\documentclass[12pt,leqno]{article}
\usepackage{amsmath,amssymb,amsbsy,amsthm}
\usepackage{latexsym,verbatim}
\usepackage[mathscr]{eucal}
\usepackage{hyperref}

\newcommand{\A}[1]{\begin{equation}\label{#1}}
\newcommand{\E}{\end{equation}}

\newcommand{\PP}{\mathbb{P}}
\newcommand{\EE}{\mathbb{E}}
\newcommand{\RR}{\mathbb{R}}

\newcommand{\la}{\langle}
\newcommand{\ra}{\rangle}

\DeclareMathOperator{\VV}{Var}
\DeclareMathOperator{\Vol}{Vol}

\newtheorem{thm}{Theorem}
\newtheorem{prop}[thm]{Proposition}
\newtheorem{lemma}[thm]{Lemma}
\newtheorem*{thm*}{Theorem}
\newtheorem*{prop*}{Proposition}
\newtheorem{cor}[thm]{Corollary}
\newtheorem*{cor*}{Corollary}
\newtheorem*{namedthm}{\theoremname}
\newcommand{\theoremname}{testing}

\theoremstyle{definition}
\newtheorem{dfn}{Definition}
\theoremstyle{remark}
\newtheorem*{rmk}{Remark}

\begin{document}

\title{Tail-sensitive Gaussian asymptotics for marginals of concentrated
measures in high dimension}

\author{Sasha Sodin}
\maketitle

\begin{abstract}

If the Euclidean norm $|\cdot|$ is strongly concentrated with respect to
a measure $\mu$, the average distribution of an average marginal of $\mu$ has
Gaussian asymptotics that captures tail behaviour.

If the marginals of $\mu$ have exponential moments, Gaussian asymptotics
for the distribution of the average marginal implies Gaussian asymptotics
for the distribution of most individual marginals.

We show applications to measures of geometric origin.

\end{abstract}

\section{Introduction}\label{S::intr}
Let $\mu$ be a probability measure on $\RR^n$; let $X \, = X_\mu$ be a random
vector distributed according to $\mu$.

We study the marginals $X^\xi = X^\xi_\mu = \la X_\mu, \, \xi \ra$ of $X_\mu$
($\xi \in S^{n-1}$); let
\[ F^\xi(t) = F^\xi_\mu(t) = \PP \{ X_\mu^\xi < t \} \]
be the distribution functions of $X_\mu^\xi$. Consider also the average marginal
$X_\mu^\text{av}$ defined by its distribution function
\[ F^\text{av}(t) = F^\text{av}_\mu(t) =
    \int_{S^{n-1}} F_\mu^\xi(t) \, d\sigma(\xi) \, \text{,} \]
where $\sigma = \sigma_{n-1}$ is the rotation-invariant
probability measure on $S^{n-1}$. If $\mu$ has no atom at the
origin, the function $F_\mu^\text{av}$ is continuously
differentiable (cf the Brehm-Voigt formul\ae{} in Section
\ref{S::aver}); denote $f_\mu^\text{av} = \left( F_\mu^\text{av}
\right)'$.

It appears that, for certain classes of measures $\mu$ on $\RR^n$, the distributions of
$X_\mu^\xi$ (for many $\xi \in S^{n-1}$) and $X_\mu^\text{av}$ are approximately Gaussian. If
$\mu = \mu_1 \otimes \mu_2 \otimes \cdots \otimes \mu_n$ is a tensor product of measures
$\mu_i$ on the real line $\RR$, this is the subject of classical limit theorems
in probability theory.

The motivation for our research comes from a different family of
measures: the (normalised) restrictions of the Lebesgue measure to
convex bodies $K \subset \RR^n$. The behaviour of the marginals of
these measures was studied recently by Anttila, Ball and
Perissinaki, Bobkov and Koldobsky, Brehm and Voigt and others
\cite{ABP,BV,BK}.

Let us state the problem more formally; denote as usual
\[  \Phi(t) = \int_{-\infty}^t \phi(s) \, ds \, \text{,}
   \quad \phi(t) = \frac{e^{-t^2/2}}{\sqrt{2\pi}} \, \text{.} \]

We wish to find sufficient conditions for proximity of distribution functions
\A{intprob} 1 - F^\text{av}(t) \approx 1 - \Phi(t)
    \, \text{,} \quad 1 - F^\xi(t) \approx 1 - \Phi(t) \, \text{,} \E
or density functions:
\A{locprob} f^\text{av}(t) \approx \phi(t)
    \, \text{,} \quad f^\xi(t) \approx \phi(t) \, \text{;} \E
we discuss the exact meaning of proximity ``$\approx$'' in the sequel. We refer to
(\ref{intprob}) as the {\em integral} problem and to (\ref{locprob}) as the {\em local}
problem.

Anttila, Ball and Perissinaki (\cite{ABP}), Brehm and Voigt (\cite{BV}), Bobkov and
Koldobsky (\cite{BK}), Romik (\cite{R}) and others proposed to study these problems
under the assumption that the Euclidean norm $|\cdot|$ is concentrated with respect
to the measure $\mu$.

These works provide a series of results, establishing
(\ref{intprob}) or (\ref{locprob}) under assumptions of this kind.
The assumptions can be verified for the geometric measures
described above (see Anttila, Ball and Perissinaki \cite{ABP}) for
some classes of bodies $K \subset \RR^n$.

However, these authors interpret ``$\approx$'' in (\ref{intprob})
and (\ref{locprob}) as proximity in $L_1$ or $L_\infty$
metrics\footnote{Recently H. Vogt \cite{V} has proved some results
concerning convergence in the $W_2^k$ Wasserstein metric.}. These
metrics fail to capture the asymptotics of the tails of the
distribution of $X^\text{av}$ beyond $t = O(\sqrt{\log n})$. We
work with a stronger notion of proximity:
\[ g \approx h \quad \text{if} \quad
    \sup_{0 \leq t \leq T} \left| \frac{g(t)}{h(t)} - 1 \right| \quad \text{is small,} \]
where $T$ may be as large as some power of $n$.

In the classical case $\mu = \mu_1 \otimes \cdots \otimes \mu_n$
this corresponds to limit theorems with moderate deviations in the
spirit of Cram\'er, Feller, Linnik et al (see Ibragimov and Linnik
\cite{IL}).

To obtain (\ref{intprob}) or (\ref{locprob}), we also assume concentration of Euclidean
norm with respect to $\mu$, but in a stronger form. That is, we reach a stronger conclusion
under stronger assumptions.

\hfill

Let us explain the results in this note. First, approach the question for average marginals
\big(the first part of (\ref{intprob}), (\ref{locprob})\big). It appears more natural to
consider ``spherical approximation'':
\[ 1 - F^\text{av}(t) \approx 1 - \Psi_n(t) \, \text{,} \quad
    f^\text{av}(t) \approx \psi_n(t) \, \text{,} \]
where
\begin{gather*}
\Psi_n(t) = \int_{-\infty}^t \psi_n(s) \, ds \, \text{,} \\
\psi_n(t) = \frac{1}{\sqrt{\pi n}} \, \frac{\Gamma(\frac{n}{2})}{\Gamma(\frac{n-1}{2})} \,
                                      \left( 1 - \frac{t^2}{n} \right)^\frac{n-3}{2} \,
                                      \mathbf{1}_{[- \sqrt{n}, \, \sqrt{n}]}(t) \, \text{.}
\end{gather*}
The geometric meaning of the distribution defined by these formul\ae{} that
justifies its name is it being the one-dimensional marginal of the uniform
probability measure on the sphere (we explain this in the proof of the
Brehm-Voigt formul\ae{} in Section \ref{S::aver}).

The following lemma shows the connection between Gaussian and spherical approximation:

\begin{lemma}\label{sph} For some constants $C, \, C_1, \, C_2 > 0$ and
some sequence $\epsilon_n \searrow 0$ the following inequalities hold\footnote{the
constant $4$ in the first inequality is written explicitly since it is sharp.}
for $0 < t < C\sqrt{n}$:
\[\begin{aligned}
(1 - \epsilon_n) \, \phi(t) \, \exp(- t^4/4n) \, &\leq \, \psi_n(t) \\
    &\leq \, (1 + \epsilon_n) \, \phi(t) \, \exp(- t^4/C_1n) \, \text{,}\\
(1 - \epsilon_n) \, (1 - \Phi(t)) \, \exp(- t^4/C_2n) \, &\leq \, 1 - \Psi_n(t) \\
    &\leq \, (1 + \epsilon_n) \, (1 - \Phi(t)) \, \exp(-t^4/C_1n) \, \text{.}
\end{aligned}\]
\end{lemma}

Informally speaking, the lemma states that Gaussian approximation for the distribution
of $X^\text{av}$ is equivalent to spherical approximation if (and only if) the variable
$t$ is small with respect to $n^{1/4}$. We prove the lemma, together with other properties
of spherical distributions, in Appendix \ref{S::sphproofs}.

Now we formulate the main result for average marginals:

\begin{thm}\label{thm_aver} Suppose for some constants $\alpha, \, \beta, \, A, \, B > 0$ we have
\begin{equation} \label{conc}
\PP \left\{ \left| \frac{|X_\mu|}{\sqrt{n}} - 1 \right| \geq u \right\}
    \, \leq \, A \, \exp \left( - B \, n^\alpha \, u^\beta \right)
\end{equation}
for $0 \leq u \leq 1$. Then
\begin{eqnarray}
&\left| \left\{1 - F^\text{av}_\mu(t)\right\} \bigg{\slash} \left\{1 - \Psi_n(t)\right\} \: - \: 1 \right|
    < C t^{2\max(\beta, \, 1)} \, n^{-\alpha} \, \text{;} \label{intmain} \\
&\left| f^\text{av}_\mu(t) \left/ \psi_n(t) \right. - 1 \right|
    < C t^{2\max(\beta, \, 1)} \, n^{-\alpha} \label{locmain}
\end{eqnarray}
for $t > 0$ s.t. $t^{2\max(\beta, \, 1)} \, n^{-\alpha}< c$; the
constants $c$, $C$ depend only on $A$, $B$, $\alpha$, $\beta$.
\end{thm}

In other words, the distribution of $X^\text{av}$ has spherical asymptotics for $t = o(n^\gamma)$, where
$\gamma = \alpha/\left( 2\max(\beta, 1) \right)$, and hence also Gaussian asymptotics for
$t = o(n^{\min(\gamma, \,1/4)})$.

We prove this theorem in Section~\ref{S::aver}.

\hfill

Then we approach the individual marginals $X^\xi$. Suppose the measure $\mu$ satisfies a
property resembling (\ref{intmain}):
\A{avappr}\begin{split}
(1 - \epsilon) \, (1 - \Phi(t))
    &\leq \int_{S^{n-1}} (1 - F^\eta(t)) \, d\sigma(\eta) \\
    &\leq (1 + \epsilon) \, (1 - \Phi(t)) \quad \text{for $0 \leq t \leq T$.} \end{split}\E
Suppose also that the measure $\mu$ has $\psi_1$ marginals:
\A{psi_1} \PP \left\{ \la X_\mu, \, \theta \ra > s \right\} \,
    \leq C \, \exp(- c s) \, \text{,}
    \quad s \in \RR^+ \, \text{,}
    \quad \theta \in S^{n-1} \, \text{.} \E

The following inequality due to Borell (see eg Giannopoulos
\cite[Section 2.1]{G} or Milman - Schechtman \cite{MS}) shows that
this property holds for an important class of measures.

\begin{dfn}
A measure $\mu$ on $\RR^n$ is called isotropic if
\A{is} \VV \la X_\mu , \, \xi \ra = 1 \quad \text{for $\xi \in S^{n-1}$ \, .} \E
\end{dfn}

\begin{dfn}
A measure $\mu$ on $\RR^n$ is called log-concave if
\A{lc} \mu \left( \frac{A + B}{2} \right) \geq \sqrt{\mu(A) \, \mu(B)}
    \quad \text{for $A, \, B \subset \RR^n$.} \E
\end{dfn}

\begin{prop*}[Borell]
Every isotropic, log-concave, even measure $\mu$ on $\RR^n$ has $\psi_1$
marginals (\ref{psi_1}).
\end{prop*}

\begin{rmk}
Actually, the isotropicity condition is too rigid, and measures satisfying a
weaker condition
\A{subis} \VV \la X_\mu , \, \xi \ra \leq C' \quad
    \text{for $\xi \in S^{n-1}$ \, .} \E
also have $\psi_1$ marginals, with constants $C$ and $c$ in (\ref{psi_1})
depending on $C'$. Such measures are called {\em ($C$-)subisotropic}.
\end{rmk}

Our aim is to show that for most $\xi \in S^{n-1}$
\A{indappr}
(1 - 10 \epsilon) \, (1 - \Phi(t)) \leq 1 - F^\xi(t)
    \leq (1 + 10 \epsilon) \, (1 - \Phi(t)) \quad \text{for $0 \leq t \leq T$;} \E
of course, the constant $10$ has no special meaning (but influences the meaning of
``most'').

This should be compared with classical results on concentration of marginal
distributions of isotropic measures.

To the extent of the author's knowledge, the earliest result of
this kind is due to Sudakov (\cite{Su}, see also von Weizs\"acker
\cite{W}). It states that if $n \geq n_0(\epsilon)$ and $\mu$ is a
general isotropic measure on $\RR^n$, then
\[ \sigma \left\{ \xi \in S^{n-1} \, \Big{|} \,
    \left\| F^\xi - F^\text{av} \right\|_1 > \epsilon \right\} \leq \epsilon \, \text{.} \]

Anttila, Ball and Perissinaki have considered isotropic measures
$\mu$ that are normalised restrictions of the Lebesgue measure to
convex bodies $K \subset \RR^n$; their work extends to general
isotropic log-concave measures. The result in \cite{ABP} states
that in this case
\[ \sigma \left\{ \xi \in S^{n-1} \,
    \Big{|} \, \left\| F^\xi - F^\text{av} \right\|_\infty > \delta \right\}
    \leq C \sqrt{n} \, \log n \, \exp \left( - c n \delta^2 \right) \, \text{.} \]

Bobkov (\cite{B}) improved both aforementioned results. In the
log-concave case he proved that for some constant $b > 0$
\[ \sigma \left\{ \xi \in S^{n-1} \,
    \Big{|} \, \sup_{t \in \RR} e^{b t} \left| F^\xi(t) - F^\text{av}(t)\right|
        > \delta \right\} \leq C \sqrt{n} \, \log n \, \exp \left( - c n \delta^2 \right)
    \, \text{.} \]
Note that the metric that appears in this inequality takes the
tails of the distributions into account. Moreover, it seems reasonable
that the term $e^{bt}$ can not be replaced by
$e^{bt^{1+\epsilon}}$ without additional assumptions.

On the other hand, the Gaussian case (\ref{avappr}) is of special
interest (see \cite{ABP,B,R,W}). The cited results allow to deduce
(\ref{indappr}) from (\ref{avappr}) only for
$T = O \left( \log^{1/2} n \right)$.

Our results show that in fact (\ref{avappr}) implies (\ref{indappr})
for $T$ as large as a certain power of $n$. Let us formulate the exact
statements.

We consider even measures with $\psi_1$ marginals.

\begin{thm}\label{intthm}
There exists $\epsilon_0 > 0$ such that if for some $\epsilon < \epsilon_0$
\[\begin{split}
(1 - \epsilon) \, (1 - \Phi(t)) &\leq \int_{S^{n-1}} (1 - F_\mu^\eta(t)) \, d\sigma(\eta) \\
    &\leq (1 + \epsilon) \, (1 - \Phi(t)) \, \text{,} \quad 0 \leq t \leq T \, \text{,}
\end{split}\]
then
\begin{multline}\label{intineq}
\sigma \left\{ \xi \in S^{n-1} \, \Big{|} \, \exists \, 0 \leq t \leq T , \,
        \left| \frac{1-F_\mu^\xi(t)}{1 - \Phi(t)} - 1\right| > 10 \epsilon \right\} \\
    \leq \frac{C T^8}{n \epsilon^4} \, \exp \left( - c \, n\epsilon^2 \, T^{-6} \right)
        \, \text{.} \hspace{2cm}
\end{multline}
\end{thm}

The constants $C$,~$c$,~$c_1$,~$\epsilon_0$,~\dots \, in this theorem, as well
as the constants in the following theorem and all other constants
in this note, depend neither on $\mu$ nor on the dimension $n$.

\begin{cor}\label{intcor}
If under assumptions of Theorem \ref{intthm}
\A{intzeta} 0 \leq T \leq \left\{ \frac{c_1 n \epsilon^2}
                              {\log n + \log \frac{1}{\epsilon} + \log \frac{1}{\zeta}}
                         \right\}^{1/6} \, \text{,} \E
then
\[ \sigma \left\{ \xi \in S^{n-1} \, \Big{|} \, \exists \, 0 \leq t \leq T , \,
        \left| \frac{1-F_\mu^\xi(t)}{1 - \Phi(t)} - 1 \right| > 10 \epsilon \right\} \leq \zeta
            \, \text{.} \]
\end{cor}

\begin{proof}[Proof of Corollary]
Substitute (\ref{intzeta}) into (\ref{intineq}). We obtain:
\[\begin{split}
\sigma \left\{ \cdots \right\}
    &\leq \frac{C}{n \epsilon^4}
         \left\{ \frac{c_1 n \epsilon^2}
                      {\log \frac{n}{\epsilon\zeta} }\right\}^{4/3}
         \exp \left\{- \frac{c}{c_1} \log \frac{n}{\epsilon\zeta}\right\} \\
    &= C c_1^{4/3} n^{1/3 - c/c_1} \epsilon^{-4/3 + c/c_1} \zeta^{c/c_1}
            \log^{-4/3} \frac{n}{\epsilon\zeta} \, \text{.}
\end{split}\]
If $c_1$ is small enough, this expression is less than $\zeta$.
\end{proof}

We also prove a local version of the theorem.  Suppose
$F_\mu^\eta$ are concave on $\RR_+$; then $f_\mu^\eta =
\left(F_\mu^\eta\right)'$ are defined a.e. and
\[ f_\mu^\text{av}(t) = \int_{S^{n-1}} f_\mu^\eta(t) \, d\sigma(\eta) \, \text{.} \]

\begin{thm}\label{locthm}
Suppose
\[ (1 - \epsilon) \, \phi(t) \leq \int_{S^{n-1}} f_\mu^\eta(t)
\, d\sigma(\eta)
    \leq (1 + \epsilon) \, \phi(t) \, \text{,} \quad 0 \leq t \leq T \, \text{.} \]
Then
\[ \sigma \left\{ \xi \in S^{n-1} \, \Big{|} \, \exists \, 0 \leq t \leq T , \,
        \left| \frac{f_\mu^\xi(t)}{\phi(t)} - 1 \right| > 10 \epsilon \right\}
    \leq \frac{C T^8}{n \epsilon^7} \,
    \exp \left( - c_1 \, n\epsilon^4 \, T^{-6} \right) \, \text{.} \]
\end{thm}

\begin{cor}
If under assumptions of Theorem \ref{locthm}
\[ 0 \leq T \leq \left\{ \frac{c_1 n \epsilon^4}
                              {\log n + \log \frac{1}{\epsilon} + \log \frac{1}{\zeta}}
                         \right\}^{1/6} \, \text{,} \]
then
\[ \sigma \left\{ \xi \in S^{n-1} \, \Big{|} \, \exists \, 0 \leq t \leq T , \,
        \left| \frac{f_\mu^\xi(t)}{\phi(t)} - 1 \right| > 10 \epsilon \right\} \leq \zeta
            \, \text{.} \]
\end{cor}

The Corollary follows from Theorem~\ref{locthm} exactly as
Corollary~\ref{intcor} follows from Theorem~\ref{intthm}.
Note that the only essential difference between the local and the
integral versions is in the dependence on $\epsilon$.

We prove the theorems in Section \ref{S::ind}. Finally, in Section
\ref{S::ex} we apply our results from Sections \ref{S::aver},
\ref{S::ind} to measures associated with convex bodies $K \subset
\RR^n$; these examples are parallel to those by Anttila, Ball and
Perissinaki \cite{ABP}.

We devote Appendix \ref{S::sphproofs} to proofs of some properties
of the spherical distribution that we use in Section \ref{S::aver}.

\hfill

{\em Acknowledgements.} I express the sincere gratitude to my supervisor
Professor Vitali Milman who introduced me to the subject, guided along
the research and encouraged to write this note. I thank Dr.~Bo'az
Klartag for many useful and encouraging discussions and for reading a
preliminary version of the text. I thank Professor Sergey Bobkov for
explaining several theorems related to concentration of marginal
distributions and for reading a preliminary version of this text.

\section{Average marginals}\label{S::aver}

We commence with explicit formul\ae{} for $1 - F^\text{av}(t)$,
$f^\text{av}(t)$, due to Brehm and Voigt (\cite{BV}, see also
Bobkov and Koldobsky \cite{BK}). Then we develop these formul\ae{}
to obtain the estimate in Proposition \ref{BV2int} (below).
Finally, we bound the integrals that appear in the estimate to
conclude the proof of Theorem \ref{thm_aver}.

Denote by $\mu^\ast$ the normalised radial projection
$$ \mu^\ast(r) = \PP \{ |X_\mu| \leq \sqrt{n} r \}
    = \mu\{B(0; \sqrt{n}r)\} \, \text{.} $$

\begin{prop*}[Brehm -- Voigt]
For any Borel probability measure $\mu$ on $\RR^n$ with
$\mu(\{0\}) = 0$, $1 - F^\text{av} \in C^1(\RR)$ and
\begin{eqnarray}
1 - F^\text{av} (t) &=&
    \int_0^\infty \left\{ 1 - \Psi_n \left( \frac{t}{r} \right) \right\} \,
    d\mu^\ast(r) \label{BV1int} \\
f^\text{av} (t) &=&
    \int_0^\infty \frac{1}{r} \, \psi_n \left( \frac{t}{r} \right) \, d\mu^\ast(r) \,
    \text{.} \label{BV1diff}
\end{eqnarray}
\end{prop*}

For completeness, we prove this proposition.

\begin{proof}[Proof of Proposition]

Proof of (\ref{BV1int}): First, let us verify the formula for $\mu = \sigma_{n-1}$.
Let us project $\sigma_{n-1}$ onto the $x$-axis; let $x_0 = \sin \theta_0$. Then
\[ \PP \left\{ x < x_0 \right\}
    = \frac{\int_{-\pi/2}^{\theta_0} \cos^{n-2} \theta d\theta}
           {\int_{-\pi/2}^{\pi/2} \cos^{n-2} \theta d\theta}
    \, \text{.} \]
Let $x = \sin \theta$, $dx = \cos \theta d\theta$; then the numerator equals
\[ \int_{-\pi/2}^{\theta_0} \cos^{n-3} \theta \, \cos \theta \, d\theta
    = \int_{-1}^{x_0} (1 - x^2)^{(n-3)/2} \, dx \, \text{.} \]
The denominator is just a constant, and the correct one, since both $\Psi_n$
and the marginal of $\sigma_{n-1}$ are probability distributions. This proves
the proposition for $\sigma_{n-1}$.

Next, let $\mu$ be a rotation-invariant measure. Then we can approximate $\mu$
by a convex combination of dilations of $\sigma_{n-1}$; these combinations satisfy
(\ref{BV1int}). Now we can pass to the limit by the dominated convergence theorem.

Finally, both sides of (\ref{BV1int}) are equal for $\mu$ and its
symmetrisation $\widetilde{\mu} = \int_{O(n)} T^*(\mu) \,
d\sigma(T)$ (here $\sigma$ is the translation-invariant measure on
the orthogonal group $O(n)$), and hence the formula extends to
arbitrary probability measures.

Proof of (\ref{BV1diff}): Apply (\ref{BV1int}) to $\mu_B = \mu(B)^{-1} \mu|_B$ for
Borel sets $B$; (\ref{BV1diff}) follows by use Fubini's theorem. To see that
$f^\text{av}$ is continuous, it suffices to check that
$$\int_0^\infty \left|\frac{d}{dt} \psi_n(t/r) \right| \, dt < \infty \, \text{.}$$
This condition can be verified by straightforward computation (cf. second statement
in Lemma \ref{sphder} in the sequel).

\end{proof}

We develop the integral formula (\ref{BV1int}) needed for the
proof of (\ref{intmain}); note that without loss of generality
$\mu$ has no atom at the origin. The computations for the local
version (\ref{locmain}) are rather similar; we prove all the
needed technical lemmata in both versions. Anyway, at the end of
the computations both questions reduce to asymptotics of the same
integral (\ref{intred}).

First, split the domain of integration in (\ref{BV1int}) into $3$ parts:
\begin{multline*}
1 - F^\text{av}(t)
    \, = \, \int_0^1 \left\{ 1 - \Psi_n \left( \frac{t}{r} \right) \right\} \, d\mu^\ast(r) \\
    \, - \, \left[ \int_1^2 + \int_2^\infty \right]
        \left\{ 1 - \Psi_n \left( \frac{t}{r} \right) \right\} \, d\left[1 - \mu^\ast(r)\right]
    \, \text{.}
\end{multline*}
Integrating by parts, we deduce:
\begin{multline*}
1 - F^\text{av}(t) \\
    \, = \, \left. \left\{ 1 - \Psi_n\left(\frac{t}{r}\right) \right\} \, \mu^\ast(r) \, \right|_0^1
    \, - \, \left. \left\{ 1 - \Psi_n\left(\frac{t}{r}\right) \right\} \,
        \left( 1 - \mu^\ast(r) \right) \, \right|_1^2 \\
    \, - \,  \int_0^1 \frac{t}{r^2} \, \psi_n \left( \frac{t}{r} \right) \, \mu^\ast(r) \, dr
    \, + \,  \int_1^2 \frac{t}{r^2} \, \psi_n \left( \frac{t}{r} \right) \, \left( 1 - \mu^\ast(r)\right) \, dr \\
    \, + \,  \int_2^\infty
        \left\{ 1 - \Psi_n \left( \frac{t}{r} \right) \right\} \, d\left[1 - \mu^\ast(r)\right]
\end{multline*}
and hence
\begin{multline*}
\left\{ 1 - F^\text{av}(t) \right\} - \left\{ 1 - \Psi_n(t) \right\}
    \, = \, - \left\{ 1 - \Psi_n(t/2) \right\} \, \left(1 - \mu^\ast(2)\right) \\
    \, - \,  \int_0^1 \frac{t}{r^2} \, \psi_n \left( \frac{t}{r} \right) \, \mu^\ast(r) \, dr
    \, + \,  \int_1^2 \frac{t}{r^2} \, \psi_n \left( \frac{t}{r} \right) \, \left( 1 - \mu^\ast(r)\right) \, dr \\
    \, + \,  \int_2^\infty
        \left\{ 1 - \Psi_n \left( \frac{t}{r} \right) \right\} \, d\left[1 - \mu^\ast(r)\right] \, \text{.}
\end{multline*}

Now we need to estimate $1 - \Psi_n(t)$. We formulate the needed property in a lemma that we prove in Appendix
\ref{S::sphproofs}.

\begin{lemma}\label{sphder}
\begin{eqnarray*}
&0 < C^{-1} \leq \frac{1 - \Psi_n(t)}{t^{-1}\psi_n(t)} \, \leq \, C \, \quad\text{for $8 t^2 < n$,}\\
&\frac{\psi_n(t)}{t^{-1}\psi_n'(t)} = 1 - t^2/n \, \text{,}
\end{eqnarray*}
where $C$ is a universal constant.
\end{lemma}

This yields the following proposition:

\begin{prop}\label{BV2int}
The following inequality holds for any Borel probability measure $\mu$ on $\RR^n$:
\[\begin{split}
&\left| \frac{1 - F^\text{av}(t)}{1 - \Psi_n(t)} - 1 \right|
    \, \leq \, \left( 1 - \mu^\ast(2) \right) \, \frac{Ct}{\psi_n(t)} \\
    &\, + \, Ct^2 \, \left\{
    \, \int_0^1 \frac{1}{r^2} \, \frac{\psi_n \left( \frac{t}{r} \right)}{\psi_n(t)} \, \mu^\ast(r) \, dr
    \, + \,  \int_1^2 \frac{1}{r^2} \, \frac{\psi_n \left( \frac{t}{r} \right)}{\psi_n(t)} \,
        \left( 1 - \mu^\ast(r)\right) \, dr \right\} \, \text{.}
\end{split}\]
\end{prop}

Now we can conclude the proof of Theorem~\ref{thm_aver}.

\begin{proof}[Proof of Theorem \ref{thm_aver}]

Apply Proposition~\ref{BV2int} and denote
\[\begin{split}
\text{TERM}_1 &= \left( 1 - \mu^\ast(2) \right) \, \frac{Ct}{\psi_n(t)} \, \text{,} \\
\text{TERM}_2 &= \int_0^1 \frac{1}{r^2} \, \frac{\psi_n \left( \frac{t}{r} \right)}{\psi_n(t)} \,
    \mu^\ast(r) \, dr \, \text{,} \\
\text{TERM}_3 &= \int_1^2 \frac{1}{r^2} \, \frac{\psi_n \left( \frac{t}{r} \right)}{\psi_n(t)} \,
    \left( 1 - \mu^\ast(r)\right) \, dr \, \text{.}
\end{split}\]

By Lemma \ref{sph} and the concentration condition (\ref{conc}) (used with $u = 1$),
\A{term1} \text{TERM}_1 \leq A' \, t \, \exp(B't^2 - B''n^\alpha) \E
with $A' = AC$, $B' = B/2$, $B'' = 2^\beta B$; this expression surely satisfies the bound (\ref{intmain}).

Introduce a new variable $u = (r - 1)$ in $\text{TERM}_3$ and use (\ref{conc}) once again. We obtain:
\begin{eqnarray*}
\text{TERM}_3 &\leq& \,
    \int_0^1 \frac{\psi_n \left(\frac{t}{1+u}\right)}{\psi_n(t)}
        \times A \exp(- B n^\alpha u^\beta) \: du \, \text{.}
\end{eqnarray*}

Now we use one more property of spherical distributions which we also prove in
Appendix \ref{S::sphproofs}.

\begin{lemma}\label{logder}
There exist constants $C_1$ and $C_2$ such that for $2t^2<n$
\[\exp(C_1 u t^2) \leq \frac{\psi_n(t)}{\psi_n((1+u)t)}\]
and for $2(1+u)^2 \, t^2 < n$
\[\frac{\psi_n(t)}{\psi_n((1+u)t)} \leq \exp(C_2 u t^2)\,\text{.}\]
\end{lemma}

By the lemma for $2 t^2 < n$
\begin{equation}\label{intred}
\text{TERM}_3 \, \leq \, A \, \int_0^1 \exp \left(C_0 t^2 u - B n^\alpha u^\beta\right)
    \, du \, \text{.}
\end{equation}
The computations in their local version would lead us to the same integral.

Now we study the integral
$$ I(K; \, L) = \int_0^1 \exp(K u - L u^\beta) \, du \, \text{,} $$
where $K$ and $L$ are large parametres, $K$ much smaller than $L$.

The exponent $E(u) = Ku - Lu^\beta$ is concave for $\beta > 1$ and convex for $\beta \leq 1$.
Let us consider these cases separately.

\paragraph{Case $1$: $\beta > 1$}\hfill

The maximum of the concave function $E(u) = Ku - Lu^\beta$ is achieved at the point
$u_0 = (K/\beta L)^{1/(\beta-1)}$ inside the domain of integration;
$E(u_0) = C_\beta (K^\beta/L)^{1/(\beta-1)}$, where
$C_\beta = [\beta^{-1/(\beta-1)} - \beta^{-\beta/(\beta-1)}]$. Let $R > 1$ be
fixed later (so that $Ru_0 \leq 1$).

First, consider the integral from $0$ to $Ru_0$.
\begin{equation}\label{add1}
\begin{split}
\int_{0}^{Ru_0} &\exp(Ku - Lu^\beta) \, du \\
    &\leq R u_0 \exp E(u_0)
    = (K/\beta L)^{1/(\beta-1)} \: R \: \exp \left\{ C_\beta (K^\beta/L)^{1/(\beta-1)} \right\} \\
    &= R \: \frac{(K^\beta/L)^{1/(\beta-1)}}{K} \:
        \frac{\exp \left\{ C_\beta (K^\beta/L)^{1/(\beta-1)} \right\}}{\beta^{1/(\beta-1)}}
\end{split}
\end{equation}

Next, for $u \geq R u_0$ we have:
\begin{multline}\label{add2}
E(u) \:\leq\:  E(u_0) \, + \, E'(Ru_0) \, (u - Ru_0) \\
    = C_\beta \, (K^\beta /L)^{1/(\beta-1)} \, - \, K \, (R^{\beta-1} - 1) \, (u - Ru_0)
\end{multline}
and hence
$$ \int_{Ru_0}^1
    \leq \frac{\exp \left(C_\beta (K^\beta/L)^{1/(\beta-1)}\right)}{K (R^{\beta-1}-1)} \, \text{.} $$

For $K^\beta/L < 1/2$ choose $R = (L/K^\beta)^{1/\beta(\beta-1)}$; then both (\ref{add1}) and
(\ref{add2}) are bounded by a constant times
    $$\frac{K^\beta/L}{K} \: \exp \left(C_\beta (K^\beta/L)^{1/(\beta-1)}\right)
        \, \leq \, C_\beta' \, \frac{K^\beta/L}{K} \, \text{.}$$

\paragraph{Case 2: $\beta \leq 1$}\hfill

For $K/L < 1/2$ the inequality $Ku \leq Lu^\beta/2$ holds in the interval $[0, \, 1]$;
hence
\begin{equation*}
\begin{split}
I(K; \, L) &\leq \int_0^1 \exp(- Lu^\beta/2) \, du
    \: = \: \beta^{-1} \, (2/L)^{1/\beta} \, \int_0^{L/2} \exp(-v) \, v^{1/\beta} \, dv \\
    &\leq \: \beta^{-1} \, (2/L)^{1/\beta} \, \Gamma(1/\beta)
    \: = \: 2^{1/\beta} \, \Gamma(\beta^{-1} + 1) \, \frac{K/L}{K} \, \text{.}
\end{split}
\end{equation*}

\vspace{10pt}

We have proved the following proposition:

\begin{prop}\label{lapl}
If $K^{\max(\beta, \, 1)}/L < 1/2$, then
$$ K \, \int_0^1 \exp(Ku - Lu^\beta)\,du \, \leq \, C \, \frac{K^{\max(\beta, \, 1)}}{L} \, \text{,} $$
where C depends only on $\beta$.
\end{prop}

Taking $K = C_0 t^2$, $L=B n^\alpha$ we arrive at the desired estimate for $\text{TERM}_3$. The integral
$\text{TERM}_2$ is even smaller, since $\psi_n(t/r)/\psi_n(t) < 1$ for $r < 1$.

\end{proof}

\section{Individual marginals}\label{S::ind}

Along the remainder of this note, we only deal with the upper
bounds in Theorems \ref{intthm} and \ref{locthm}. The same
technique works also for lower bounds. Note that these bounds
do not depend on each other: the left side inequality in
(\ref{avappr}) implies the left side inequality in
(\ref{indappr}), and similarly for the right side inequalities.

Also, all the measures $\mu$ in this section are assumed
even with $\psi_1$ marginals; we reiterate that all the
constants do not depend on $\mu$ nor on the dimension $n$.

Let us explain the idea of the proof (of the integral theorem).
Let $A$ be the set of directions $\eta \in S^{n-1}$ such that $1 -
F^\eta(t - s)$ is not too large. Markov's inequality combined with
the bound (\ref{avappr}) for the average marginal shows that the
measure of $A$ is not too small.

Now use the triangle inequality in the following form:
\A{triang}\begin{split}
1 - F^\eta(t + s) &- \PP \left\{ X^\xi - X^\eta > s \right\}
        \leq 1 - F^\xi(t) \\
    &\leq 1 - F^\eta(t - s)
        + \PP \left\{ X^\xi - X^\eta > s \right\} \, \text{;}
\end{split}\E
we need the right side for the upper bounds.

Consider directions $\xi$ in the $\delta$-extension of $A$
\[ \left\{A\right\}_\delta = \left\{ \xi \in S^{n-1} \, \big{|} \,
    \exists \, \eta \in A, \, | \xi - \eta | \leq \delta \right\} \]
For such $\xi$, the term $1 - F^\eta(t - s)$ is not too large; the
term $\PP \left\{ X^\xi - X^\eta > s \right\}$ can be bounded in
terms of $\delta$ using the $\psi_1$ condition (\ref{psi_1}).

Finally, we use the spherical isoperimetric inequality to show
that $\left\{ A \right\}_\delta$ covers most of the sphere.

\hfill

Now we pass to rigorous exposition of the idea explained above. Define
the set of "good directions"
\[ A \left(t; \, \epsilon\right) =
    \left\{ \eta \in S^{n-1} \, \mid \,
        (1 - F^\eta(t) \leq (1 - \Phi(t)) \,
        (1 + \epsilon) \right\} \, \text{.} \]

Our first aim is to prove the following proposition:

\begin{prop}\label{ext}
There exist $c, \epsilon_0 > 0$ (that depend neither on
$\mu$ nor on $n$) such that for every $t > 1$ there exists $t' < t$ satisfying
\[ \left\{ A \left( t'; \, \epsilon \right)\right\}_{c \epsilon t^{-3}}
\subset A_{t;  \, 4 \epsilon}
\quad \text{for $0 \leq \epsilon \leq \epsilon_0$.} \]
\end{prop}

\begin{rmk}
Note that all the results for $0 \leq t \leq 1$ follow from the known
results (for example, \cite{ABP}), and hence we may restrict ourselves
to $1 \leq t$ along all the proofs.
\end{rmk}

\begin{proof}[Proof of Proposition \ref{ext}]

Suppose $1 - F^\eta(t - s) \leq (1 - \Phi(t - s)) \, (1 + \epsilon)$.
Combining (\ref{triang}) with the $\psi_1$ condition (\ref{psi_1}), we deduce:
\begin{equation}\label{tr1}
1 - F^\xi(t) \leq (1 - \Phi(t - s)) \, (1 + \epsilon) + C \,
\exp(- c s \delta^{-1} ) \, \text{,}
\end{equation}
where $\delta = | \xi - \eta |$.

Now we need to use properties of the Gaussian distribution that
are summarised in the following elementary lemma (cf Lemmata
\ref{sphder} and \ref{logder}):
\begin{lemma}\label{normal}
The following inequalities hold:
\begin{eqnarray}
1 - \Phi(t - s) &\leq& (1 - \Phi(t)) \, \exp(st) \, \text{;} \\
1 - \Phi(t) &\leq& C \, t^{-1} \, \exp(-t^2/2) \, \text{.}
\end{eqnarray}
\end{lemma}

Substituting these inequalities into (\ref{tr1}), we obtain:
\begin{equation}\label{mainineq}
\frac{1 - F^\xi(t)}{1 - \Phi(t)}
    \leq (1 + \epsilon) \, e^{st} + C_1 t \,
    \exp \left[ t^2/2 - c \, s \, \delta^{-1} \right]  \, \text{.}
\end{equation}

This inequality holds for any $s > 0$, and $s$ does not appear on
its left side. To conclude the proof, we optimise over $s$ in a
rather standard way. Denote
\[ a(s) = (1 + \epsilon) \, e^{st}
    + C_1 t \exp\left[ t^2/2 - c \, s \, \delta^{-1} \right] \, \text{.} \]
Then
\[ a'(s) = t \, \left\{ (1 + \epsilon) \, e^{st}
    - C_1 c  \, \delta^{-1} \, \exp\left[ t^2/2 - c \, s \, \delta^{-1} \right] \right\} \]
and hence the minimum is obtained at $s_0$ such that
\[ (1 + \epsilon) \, e^{s_0t} =
    C_1 c  \, \delta^{-1} \, \exp\left[ t^2/2 - c \, s_0 \, \delta^{-1} \right] \, \text{,} \]
or:
\begin{equation}\label{s0}
e^{s_0t}
    = \left( \frac{C_1 c \, e^{t^2/2}}{(1 + \epsilon) \, \delta} \right)^
        {\left[{1 + c \, \delta^{-1} \, t^{-1}}\right]^{-1}} \, \text{.}
\end{equation}

Hereby
\[ a(s_0) = (1 + \epsilon) \, \left(1 + \frac{t \delta}{c}\right) \,
    \left( \frac{C_1 c \, e^{t^2/2}}{(1 + \epsilon) \, \delta} \right)^
        {\left[{1 + c \, \delta^{-1} \, t^{-1}}\right]^{-1}} \, \text{.} \]

Extracting logarithms, we see that
\[ \log a(s_0) \leq \epsilon + \frac{t\delta}{c}
    + \frac{C_2}{1 + c \delta^{-1} t^{-1}}
    + \frac{t^2/2}{1 + c \delta^{-1} t^{-1}} \, \text{.} \]

If $\delta < c_3 \epsilon t^{-3}$, the fourth term is bounded by
$\epsilon$. If $t \geq 1$, the preceding two terms are ignorable
(and in particular their sum is bounded by $2\epsilon -
8\epsilon^2$). Finally, exploiting the inequality $\exp(u - u^2/2)
\leq 1 + u$ we deduce that $a(s_0) \leq 4 \epsilon$.

Hence $t' = t - s_0$ satisfies the requirements of the proposition.

\end{proof}

We also outline the proof of a local version of Proposition \ref{ext}. Define
\[ B \left(t; \, \epsilon\right) =
    \left\{ \eta \in S^{n-1} \, \mid \,
        f^\eta(t) \leq \phi(t) \, (1 + \epsilon) \right\} \, \text{.} \]

\begin{prop}\label{locext}
Suppose $F_\mu^\eta$ are concave on $\RR_+$. Then there exist
constants $c, \epsilon_0 > 0$ (that depend neither on $\mu$ nor on
$n$) such that for every $t > 1$ there exists $t' < t$ satisfying
\[ \left\{ B \left( t'; \, \epsilon \right)\right\}_{c \epsilon^2 t^{-3}}
    \subset B \left(t;  \, 4 \epsilon\right)
    \quad \text{for $0 \leq \epsilon \leq \epsilon_0$.} \]
\end{prop}

\begin{proof}[Sketch of proof]
Choose two small parametres, $1 \gg h \gg s > 0$. By the intermediate value theorem
\[\begin{split}
h f^\xi(t) &\leq F^\xi(t) - F^\xi(t - h) \\
    &\leq F^\eta(t + s) - F^\eta(t - h - s)
        + 2 \PP \left\{ \PP \la X, \, \xi - \eta \ra > s \right\} \\
    &\leq \left[ h + 2s \right] f^\eta(t - h - s)
        + 2 C \, \exp\left( - c \, \delta^{-1} s \right) \, \text{;}
\end{split}\]
therefore if $f^\eta(t - h - s) < (1 + \epsilon) \, \phi(t - h - s)$,
\[ \frac{f^\xi(t)}{\phi(t)}
    \leq \left[ 1 + 2 sh^{-1} \right] \, (1 + \epsilon) \, \exp\left( t (h + s) \right)
        + 2 C \exp\left( t^2/2 - c \delta^{-1} s \right) \, \text{.} \]

Take $s = C_1 \epsilon^2 t^{-1}$, $h = C_2 \epsilon t^{-1}$. For
appropriate choice of the constants $C_1$, $C_2$ we deduce:
$f^\xi(t)/\phi(t) \leq 4 \epsilon$.
\end{proof}

Now we are ready to prove Theorems \ref{intthm} and \ref{locthm}. The proofs of
these theorems are rather similar; let us prove for example the (upper bound in)
Theorem~\ref{intthm}.

\begin{proof}[Proof of Theorem \ref{intthm}]

First, apply Markov's inequality to the right side of
(\ref{avappr}). We deduce:
\[ \sigma \left\{ \eta \, \mid \, (1 - F^\eta(t) \leq (1 - \Phi(t)) \,
        (1 + 2 \epsilon) \right\}
    \geq \frac{\epsilon}{1 + 2 \epsilon} \geq \epsilon/2 \, \text{.} \]
Surely, this inequality also holds with $t'$ instead of $t$. Now
we need to transform Proposition \ref{ext} into a lower bound on
the measure of $A \left( t; \, 8 \epsilon \right)$.

Let $1 \leq t_1 \leq \cdots \leq t_I = T$ be an increasing sequence of points such that
\[ \sigma \left\{ \xi \, \mid \,
        1 - F^\xi(t_i) \geq (1 - \Phi(t)) \, (1 + 8 \epsilon) \right\}
    \leq \zeta_i \, \text{,} \quad 1 \leq i \leq I \, \text{.} \]
Then
\[ \sigma \left\{ \xi \, \mid \,
        \exists \, 1 \leq i \leq I, \: 1 - F^\xi(t_i)
            \geq (1 - \Phi(t)) \, (1 + 8 \epsilon) \right\}
    \leq \sum_{i=1}^I \zeta_i \, \text{.} \]

The function $F^\xi$ is monotone for every $\xi$; hence for $t_i \leq t \leq t_{i+1}$
we have $1 - F^\xi(t) \leq 1 - F^\xi(t_i)$. Applying Lemma \ref{normal}, we conclude:
\begin{multline*}
\sigma \left\{ \xi \, \mid \,
        \exists \, 1 \leq t_i \leq t \leq t_{i+1} \leq t_I, \right. \\
    \left. 1 - F^\xi(t) \geq \exp \left(t_{i+1}(t_{i+1}-t_i)\right) (1 - \Phi(t)) \,
        (1 + 8 \epsilon) \right\}
    \leq \sum_{i=1}^I \zeta_i \, \text{.}
\end{multline*}

Choose $t_i = \sqrt{C\epsilon i}$ with $C$ such that
$\exp \left( (t_{i+1} - t_i) \, t_{i+1} \right) \leq \epsilon$. Then
\[ \sigma \left\{ \xi \, \mid \,
        \exists \, 1 \leq t \leq T, \right. \\
    \left. 1 - F^\xi(t) \geq \exp (1 - \Phi(t)) \, (1 + 10 \epsilon) \right\}
    \leq \sum_{i=1}^I \zeta_i \, \text{.} \]

Now we use the concentration inequality on the sphere in the following form:

\begin{prop*}[Concentration on the sphere]
For $A \subset S^{n-1}$
\A{sphere_conc} \sigma (A) \,
    \left[ 1 - \sigma \left( \left\{ A \right\}_\gamma \right) \right]
    \leq \exp \left( - (n-1) \,  \gamma^2/4 \right) \, \text{.} \E
\end{prop*}

This is a standard corollary of the isoperimetric inequality on
the sphere due to P. L\'evy that can be verified applying the
concentration inequality as in Milman - Schechtman \cite{MS} to
the function $x \mapsto \inf_{y \in A} d(x, \, y)$.

Proposition \ref{ext} combined with the concentration inequality yields
\[ \zeta_i = \frac{C}{\epsilon} \,
    \exp \left[ - c_1 n \epsilon^{-1} i^{-3} \right] \, \text{;} \]
hence
\begin{align*}
\sigma \, \biggl\{ \xi \: &\bigg{|} \:
        \exists \, 1 \leq t \leq T, \,
            1 - F^\xi(t) \geq \exp \left(t_{i+1}(t_{i+1}-t_i)\right) (1 - \Phi(t)) \,
                (1 + 10 \epsilon) \biggr\} \\
    &\leq \sum_{i=1}^I \frac{C}{\epsilon}
        \exp \left[ - c_1 n \epsilon^{-1} i^{-3} \right]
    \leq \frac{C}{\epsilon}
        \int_0^I \exp \left( - c_1 n \epsilon^{-1} x^{-3} \right) \, dx \, \text{;}
\end{align*}
the second inequality is justified since the function
$i \mapsto \exp \left( - c_1 n \epsilon^{-1}i^{-3} \right)$
is monotone decreasing.

Continuing the inequality and replacing $y^{-4/3}$ with its value at the left end of
the integration domain, we obtain:
\[ \cdots \leq \frac{C_1}{\epsilon^{4/3}}
        \int_{\epsilon^2 T^{-6}}^\infty \exp(- c_1 n y) \, y^{-4/3} \, dy
    \leq C_2 \epsilon^{-4} T^8 n^{-1} \,
        \exp \left( - c_1 n \epsilon^2 T^{-6} \right)  \, \text{.} \]

\end{proof}

We conclude with a remark.

\begin{rmk}
One can generalise the conclusion of Theorems \ref{intthm} and
\ref{locthm} to measures $\mu$ satisfying the $\psi_\alpha$
property \A{psi} \PP \left\{ \la X, \, \theta \ra > s \right\} \,
    \leq C \, \exp(- c s^\alpha) \, \text{,} \quad s \in \RR^+ \E
for some $0 < \alpha \leq 2$. In this case we use (\ref{psi}) instead of (\ref{psi_1})
in the proofs of Propositions \ref{ext}, \ref{locext}. This yields $t^{-1 - 2/\alpha}$
instead of $t^{-3}$ in these propositions, leading to exponent $2 + 4\alpha^{-1}$ instead
of $6$ in the theorems.
\end{rmk}

\section{Examples}\label{S::ex}

Let us show some examples where our results apply. Our examples have geometric
motivation, hence we recall some geometric notions.

Let $K \subset \RR^n$ be a symmetric convex body; denote its boundary by $\partial K$.
Define three measures associated with $K$, called the {\em volume measure}, the
{\em surface measure} and the {\em cone measure} and denoted by $\mathscr{V}_K$,
$\mathscr{S}_K$ and $\mathscr{C}_K$ respectively:
\[\begin{split}
\mathscr{V}_K (A) &= \frac{\Vol A \cap K}{\Vol K} \, \text{;} \\
\mathscr{S}_K (A) &= \lim_{\epsilon \to + 0}
    \frac{\Vol \left\{ A \cap \partial K \right\}_\epsilon }
         {\Vol \left\{ \partial K \right\}_\epsilon} \, \text{;} \\
\mathscr{C}_K (A) &= \frac{\Vol \left\{ x \in K \, \big{|} \, x / \|x\|_K \in A \right\} }
                          {\Vol K} \, \text{.}
\end{split}\]

Here subscript denotes metric extension in $\RR^n$:
$$ \left\{ A \right\}_\epsilon =
    \left\{ a \in \RR^n \, \big{|} \, \exists \, x \in A, | x - a |
        \leq \epsilon \right\} \, \text{.} $$

\begin{rmk}
The Brunn--Minkowski inequality (see \cite{G,MS}) shows that the measure
$\mathscr{V}_K$ is log-concave for any convex body $K$.
\end{rmk}

\begin{dfn}
The body $K$ is called isotropic (subisotropic) if the measure
$\mathscr{V}_K$ is isotropic (subisotropic).
\end{dfn}

We are mainly interested in the volume measure $\mathscr{V}_K$; however, sometimes
it is easier to verify the concentration condition (\ref{conc}) for $\mathscr{S}_K$
or $\mathscr{C}_K$. As well known, the difference is insignificant:

\begin{prop}\label{sc2v}
Suppose one of the following two inequalities holds:
\begin{align}
\mathscr{S}_K \left\{ \bigl| |X| - 1 \bigr| \geq u \right\}
    \leq A \exp \left( - B n^\alpha u^\beta \right) \, &\text{,}
        &0 \leq u \leq 1& \label{surf} \\
\mathscr{C}_K \left\{ \bigl| |X| - 1 \bigr| \geq u \right\}
    \leq A \exp \left( - B n^\alpha u^\beta \right) \, &\text{,}
        &0 \leq u \leq 1&  \, \text{.} \label{cone}
\end{align}
Then
\A{vol} \mathscr{V}_K \left\{ \bigl| |X| - 1 \bigr| \geq u \right\}
    \leq A' \exp \left( - B' n^{\min(\alpha, \, 1)} u^{\max(\beta, \, 1)} \right) \, \text{,}
    \quad 0 \leq u \leq 1  \, \text{,} \E
where $A'$, $B'$ depend only on $A$, $B$, $\alpha$, $\beta$.
\end{prop}

\begin{proof}

Let $X$ be distributed according to $\mathscr{V}_K$; then $X/\|X\|_K$ is
distributed according to $\mathscr{C}_K$ and
$\PP \left\{ \|X\|_K \leq r \right\} = r^n$ for $0 \leq r \leq 1$.

Similarly, if $Y$ is distributed according to $\mathscr{S}_K$ and $R$ is a (scalar)
random variable that does not depend on $Y$ such that
\[ \PP \left\{ R\leq r \right\} = r^n \quad \text{for $0 \leq r \leq 1$,} \]
then $RY$ is distributed according to $\mathscr{V}_K$.

Therefore
\begin{equation}\label{surf1}\begin{split}
\mathscr{V}_K \left\{ |x| < 1 - u \right\}
    &\leq \mathscr{V}_K \left\{ |x| < (1 - u/2)^2 \right\} \\
    &\leq \mathscr{S}_K \left\{ |x| < (1 - u/2) \right\} + (1 - u/2)^n \\
    &\leq \mathscr{S}_K \left\{ |x| < (1 - u/2) \right\} + \exp(-nu/2)
\end{split}\end{equation}
and also
\A{cone1} \mathscr{V}_K \left\{ |x| < 1 - u \right\}
    \leq \mathscr{C}_K \left\{ |x| < (1 - u/2) \right\} + \exp(-nu/2) \, \text{.} \E

On the other hand,
\begin{align}
\mathscr{V}_K \left\{ |x| > 1 + u \right\} &\leq \mathscr{S}_K \left\{ |x| > (1 + u) \right\}
    \, \text{,} \label{surf2}\\
\mathscr{V}_K \left\{ |x| > 1 + u \right\} &\leq \mathscr{C}_K \left\{ |x| > (1 + u) \right\}
    \, \text{.} \label{cone2}
\end{align}

Combining (\ref{surf}) with (\ref{surf1}) and (\ref{surf2}) or (\ref{cone})
with (\ref{cone1}) and (\ref{cone2}), we arrive at (\ref{vol}).

\end{proof}

We also note that sometimes for $K$ in natural normalisation we
get \A{conc1} \mathscr{V}_K \left\{ x \in \RR^n \, \big{|} \,
    \left| \frac{|x|}{C_K\sqrt{n}} - 1 \right| \geq u \right\}
    \leq A' \exp \left( - B' n^{\alpha} u^{\beta} \right) \E
for $0 \leq u \leq 1$, instead of (\ref{conc}). Then we obtain spherical
asymptotics for the distribution of $X_{\mathscr{V}_K}^\text{av}/C_K$
instead of $X_{\mathscr{V}_K}^\text{av}$.

\subsection{The \(l_p\) unit balls}\label{s::lp}

The result of this subsection is

\begin{cor}\label{lp}
For $1 \leq p \leq \infty$ the average marginal of $\mathscr{V}_{B_p^n}$ has
Gaussian asymptotics for $t = o \left (n^{1/4} \right)$. Almost all marginals
of $\mathscr{V}_{B_p^n}$ have Gaussian asymptotics for
\[ t = o \left( \left[ \frac{n}{\log n}
    \right]^{1/\{2 + 4/\min(p, \, 2)\}} \right) \, \text{.} \]
\end{cor}

\begin{rmk}[Rigorous meaning of Corollary \ref{lp}]\hfill
\begin{enumerate}
\item\label{lp::1} Writing ``Gaussian asymptotics of a random variable~$X$'',
we really mean Gaussian asymptotics for $X/C$ for some $C > 0$. The power
$1/\{2 + 4/\min(p, \, 2)\}$ is between $1/6$ and $1/4$.
\item\label{lp::2} It seems natural to take $C = \sqrt{\VV X}$; however,
strictly speaking, this can not be done under general assumptions
(for a general body $K$). In the special case of $B_p^n$ one can
combine the inequality (\ref{allp}) (below) with an inequality for
$u \geq 1$ and then use the methods described in Milman -
Schechtman \cite[Appendix V]{MS} to show that one can take $C_p =
\sqrt{\VV X^\xi_{\mathscr{V}_{B_p^n}}}$ in (\ref{allp}) without
loss of generality. Here $\xi$ is of no importance, since $\VV
X^\xi_{\mathscr{V}_{B_p^n}}$ does not depend on $\xi \in S^{n-1}$.
We pay no further attention to these issues.
\item\label{lp::3} The rigorous meaning of the expression ``almost all marginals''
is as in Theorems \ref{intthm} and \ref{locthm}.
\end{enumerate}
\end{rmk}

To prove the first part of this corollary, we verify (\ref{conc}) (or, rather, (\ref{conc1}))
for $\mathscr{C}_K$, where $K = B_p^n$ is the $l_p^n$ unit ball.

For $2 \leq p < \infty$, a reasonable estimate can be obtained
using the representation of $\mathscr{C}_{B_p^n}$ found by
Schechtman and Zinn and independently by Rachev and R\"uschendorf
(\cite{SZ1,RR}; see Barthe, Gu\'edon, Mendelson and Naor
\cite{BGMN} for an extension to $\mathscr{V}_{B_p^n}$).

\begin{thm*}[Schechtman -- Zinn, Rachev -- R\"uschendorf]
Let $g_1$,~\dots,~$g_n$ be independent identically distributed random variables with density
\[ \left( 2 \Gamma(1 + p^{-1}) \right)^{-1} \exp(- |t|^p) \, \text{.} \]
Denote $G = (g_1, \, \dots, \, g_n)$ and consider the random vector $V = G/\|G\|_p$.
Then $V$ is distributed according to $\mathscr{C}_{B_p^n}$.
\end{thm*}

\begin{cor*} For $2 \leq p < \infty$ the inequality
\A{bigp_cone} \mathscr{C}_{B_p^n}
    \left\{ \left| \|V\|_2 \left/ \frac{(\EE g^2)^{1/2} \, n^{1/2}}{(\EE g^p)^{1/p} \, n^{1/p}} - 1 \right. \right|
        > u \right\} \leq A \, \exp (- B n u^2) \, \text{.} \E
holds for $0 \leq u \leq 1$.
\end{cor*}

\begin{proof}[Proof of Corollary]

The inequality
\[ (1 + u/4) \leq (1 - u/2)(1 + u) \, \quad 0 \leq u \leq 1\]
implies
\[ \PP \Biggl\{ \|V\|_2 > \frac{(\EE g^2)^{1/2} \, n^{1/2}}{(\EE g^p)^{1/p} \, n^{1/p}} \: (1 + u) \Biggr\}
    \leq \PP \Biggl\{ \frac{\|G\|_2}{\|G\|_p}
        > \frac{(\EE g^2)^{1/2} \, n^{1/2} \, (1 + u/4)}
            {(\EE g^p)^{1/p} \, n^{1/p} \, (1 - u/2)} \Biggr\} \, \text{.} \]
Then,
\begin{align*}
\PP \left\{ \|G\|_2 > (\EE g^2)^{1/2} \, n^{1/2} \, (1 + u/4) \right\}
    &\leq \PP \left\{ \sum_i (g_i^2 - \EE g_i^2) > \frac{\EE g_i^2}{2} \, n u \right\} \, \text{,} \\
\PP \left\{ \|G\|_p < (\EE g^p)^{1/p} \, n^{1/p} \, (1 - u/2) \right\}
    &\leq \PP \left\{ \sum_i (g_i^p - \EE g_i^p) < - \frac{\EE g_i^p}{2} \, n u \right\} \, \text{.}
\end{align*}

Now we need an inequality due to S.~N.~Bernstein (\cite{Be}; see
Bourgain, Lindenstrauss and Milman \cite{BLM} for available
reference).

\begin{thm*}[S.~Bernstein]
Suppose $h_1$, \dots, $h_n$ are independent random variables such that
\A{psi1} \EE h_j = 0 \, \text{;} \quad
\EE \exp \left( h_j / C \right) \leq 2 \, \text{.} \E
Then
\[ \PP \left\{ h_1 + \cdots + h_n > \epsilon n \right\}
    \leq \exp \left( - \frac{\epsilon^2 n}{16 C^2} \right) \, \text{,}
    \quad 0 \leq \epsilon \leq c \sqrt{n} \, \text{.} \]
\end{thm*}

It is easy to verify that $g_i^2 - \EE g_i^2$ and $- g_i^p + \EE g_i^p$ satisfy
(\ref{psi1}) (with some constant $C$); this yields
\[ \PP \left\{ \|V\|_2 > \frac{(\EE g^2)^{1/2} \, n^{1/2}}
                              {(\EE g^p)^{1/p} \, n^{1/p}} \: (1 + u)\right\}
    \leq \frac{A}{2} \, \exp (- B n u^2) \]
for some constants $A$ and $B$. A bound for the probability of negative deviation
can be obtained in a similar way.

The estimate (\ref{bigp_cone}) follows.

\end{proof}

For $1 \leq p \leq 2$, we use the following theorem due to Schechtman and Zinn (\cite{SZ2}):

\begin{thm*}[Schechtman -- Zinn]
There exist positive constants $C$, $c$ such that if $1 \leq p \leq 2$ and $f: \partial B_p^n \to \RR$
satisfying
\A{lip} |f(x) - f(y)| \leq |x - y| \quad \text{for all $x, y \in B_p^n$} \E
then, for all $u > 0,$
\A{sz2} \mathscr{C}_{B_p^n} \left\{ x \, \big{|} \,
    \left| f(x) - \int f \, d\mathscr{C}_{B_p^n} \right| > u \right\} \leq C \exp \left( - c n u^p \right)
    \, \text{.} \E
\end{thm*}

The condition (\ref{lip}) surely holds for $f = |\cdot|$; hence
$|\cdot|$ satisfies (\ref{sz2}). For correct normalisation recall
that \A{ip} c_1 n^{1/2 - 1/p} \leq \int |x| \,
d\mathscr{C}_{B_p^n}(x) \leq c_2 n^{1/2 - 1/p} \, \text{;} \E
hence for $1 \leq p \leq 2$ \A{smallp_cone} \mathscr{C}_{B_p^n}
\left\{ x \, \big{|} \,
    \left| |x| \left/ \int |x| \, d\mathscr{C}_{B_p^n}(x) \right. - 1 \right|
        > u \right\} \leq C \exp \left( - c n^{p/2} u^p \right)
    \, \text{.} \E

\begin{cor*}[Concentration of $|\cdot|$ with respect to $\mathscr{V}_{B_p^n}$]
For $1 \leq p \leq \infty$ there exist $A_p, B_p, C_p > 0$ such that the inequality
\A{allp} \mathscr{V}_{B_p^n} \left\{ x \, \big{|} \,
    \biggl| |x| \left/ C_p \right. - 1 \biggr| > u \right\}
        \leq A \, \exp \left( - B n^{\min \left(p , \, 2 \right)/2} \,
            u^{\min \left(p , \, 2 \right)} \right) \E
holds for $0 \leq u \leq 1$.
\end{cor*}

\begin{proof}[Proof of Corollary]

For $1 \leq p \leq 2$ combine (\ref{smallp_cone}) with Proposition \ref{sc2v}.
For $2 \leq p < \infty$ combine (\ref{bigp_cone}) with Proposition \ref{sc2v}.

For $p = \infty$ the coordinates of a random vector $X = (X_1, \cdots, X_n)$
distributed according to $\mathscr{V}_{B_p^n}$ are independent; hence Bernstein's
inequality for $X_i^2 - \EE X_i^2$ yields the result.

\end{proof}

Now we can prove Corollary \ref{lp}:

\begin{proof}[Proof of Corollary \ref{lp}]

For the first statement, apply Theorem \ref{thm_aver} using (\ref{allp}) for
$\widetilde{\mathscr{V}}_{B_p^n}$,
\[ \widetilde{\mathscr{V}}_{B_p^n}(A) = \mathscr{V}_{B_p^n}(C_p A) \, \text{.} \]

For the second statement, note that the measure $\widetilde{\mathscr{V}}_{B_p^n}$
satisfies the $\psi_\alpha$ condition (\ref{psi}) with $\alpha = \min(p, \, 2)$.
Applying Theorems~\ref{intthm} and \ref{locthm} (combined with the concluding
remark in Section \ref{S::ind}) we obtain the result.

\end{proof}

\begin{rmk}
Note that Corollary \ref{lp} does not capture the change of asymptotic behaviour
that probably occurs around $t = n^{1/4}$. This is because the bound (\ref{allp})
is not sharp.
\end{rmk}

To emphasise this point, let us consider the case $p = 2$. The
surface measure $\mathscr{C}_{B_2^n}$ surely satisfies
(\ref{surf}) with any $\alpha, \, \beta > 0$; hence
$\mathscr{V}_{B_2^n}$ satisfies (\ref{vol}) with $\alpha = \beta =
1$ (as we could have also verified by direct computation).
Applying Theorem~\ref{thm_aver}, we obtain spherical asymptotics
for $t = o \left( n^{1/2} \right)$; in particular, we capture the
breakdown of Gaussian asymptotics around $t \approx n^{1/4}$
(recall Lemma \ref{sph}).

\begin{rmk}

In fact, the bound (\ref{sz2}) for the concentration of Euclidean norm with respect to
$\mathscr{C}_{B_p^n}$ is not sharp. Schechtman and Zinn proved a better bound for $p = 1$
(for $f = |\cdot|$) in the same paper \cite{SZ2}, and Naor (\cite{N}) extended their
results to all $1 \leq p \leq 2$.

Unfortunately, these bounds do not suffice to improve the result in Corollary \ref{lp}.
On the other hand, the bounds in \cite{SZ2,N} were proved exact only on part of the range
of $u$; this makes it tempting to conjecture spherical approximation for
$t = o \left( n^{\left(p^2 - p + 2\right)/8} \right)$, $1 \leq p \leq 2$.

This would be an improvement of Corollary \ref{lp} for all $1 < p < 2$; in particular,
we would be able to capture the breakdown of Gaussian asymptotics around $t = n^{1/4}$
for all these $p$.

\end{rmk}
\hfill

Now we compare these results to limit theorems with moderate
deviations for independents random variables. This allows to
analyse the sharpness of the result in Corollary \ref{lp} for the
common case $p = \infty$.

The following more general statement follows from our results:

\begin{thm}\label{prodthm}
Let $\mu = \mu_1 \otimes \cdots \otimes \mu_n$ be a tensor product of $1$-dimensional
even measures that satisfy
\A{psi2} \int_{-\infty}^{+\infty}
    \exp \left(x^2/C^2\right) \, d\mu_i(x) \leq 2 \, \text{.} \E
Then the average marginal of $\mu$ has Gaussian asymptotics for
$t = o\left(n^{1/4}\right)$. Almost all marginals of $\mu$ have Gaussian
asymptotics for $t = o\left((n/\log n)^{1/4}\right)$.
\end{thm}

\begin{rmk} The remarks \ref{lp::1}.\ and \ref{lp::3}.\ after Corollary~\ref{lp} are
still valid. On the other hand, the variance of the approximating Gaussian variable
is ``correct'' in this case.
\end{rmk}

The classical limit theorems with moderate deviations (see Feller
\cite[Chapter XV]{F} or Ibragimov - Linnik \cite{IL} for a more
general treatment) assume a weaker assumption \A{psi1'}
\int_{-\infty}^{+\infty} \exp \left(x/C\right) \, d\mu_i(x) \leq 2
\E and establish Gaussian asymptotics of $1-F^\xi(t)$ and
$f^\xi(t)$ for $t = o \left( \|\xi\|_\infty^{-1/2} \right)$; these
results are sharp. The $l_\infty$ norm of a typical vector $\xi
\in S^{n-1}$ is of order $\sqrt\frac{\log n}{n}$; hence the
asymptotics for random marginals in Theorem \ref{prodthm} is valid
for $t = o \left( \sqrt[4]{\frac{n}{\log n}} \right)$ and our
results are sharp. In particular, this is true for $p = \infty$ in
Corollary \ref{lp}.

\subsection{Uniformly convex bodies contained in small Euclidean balls}\label{ex::uc}

Let $K \subset \RR^n$ be a convex body; define the modulus of convexity
\[ \delta_K(\epsilon)
    = \min \left\{ 1 - \frac{\|x + y \|_K}{2} \, \Big{|} \,
            \|x\|_K = \|y\|_K = 1, \, \|x - y\|_K \geq \epsilon \right\}
\, \text{.} \]

The following concentration property was proved by Gromov and
Milman (\cite{GM}, see also Arias~de~Reyna, Ball and Villa
\cite{ABV}):
\begin{thm}[Gromov -- Milman]
If $A\subset K$ has positive measure, and $d_K(x,A)$ is the distance from $x$ to $A$
(measured in the norm with unit ball $K$), then
\A{GM} \mathscr{V}_K \left\{ x \, \big{|} \,
    d_K(x, \, A)>\epsilon \right\} < \frac{e^{-2n\delta_K (\epsilon )}}{\mathscr{V}_K(A)}
    \, \text{.} \E
\end{thm}

\begin{cor}

Suppose an isotropic body $K$ satisfies
\[ K \subset C n^\nu B_2^n, \quad \delta_K(\epsilon) \geq c \epsilon^\mu
    \quad \text{and} \quad \mathscr{V}_K \left\{ |x| < m \right\} > 1/2 \]
for some constants $C$, $c$, $m$. Then
\begin{enumerate}
\item the average marginal of $\mathscr{V}_K$ has spherical asymptotics for
\[ t = o\left(n^{(1/2 + \mu^{-1} - \nu)/2}\right) \, \text{;} \]
\item almost all marginals of $\mathscr{V}_K$ have Gaussian asymptotics for
\[ t = o\left(\left[\frac{ n }{ \log n }\right]^{\min \left( 1/6, \,
                (1/2 + \mu^{-1} - \nu)/2 \right)} \right) \, \text{.} \]
\end{enumerate}

\end{cor}

\begin{proof}[Proof of Corollary]
Following \cite{ABP} we show that (\ref{GM}) implies concentration
of the Euclidean norm. Really, one can estimate the probability of
deviation from the median $\mathcal{M}$:
\begin{multline}\label{uc1}
\mathscr{V}_K \left\{ |x| \leq \mathcal{M} - \epsilon \right\}
    \leq \mathscr{V}_K \left\{ d_2 \left( x, \left\{ y \leq \mathcal{M} \right\} \right)
        > \epsilon \right\} \\
    \leq \mathscr{V}_K \left\{ d_K \left( x, \left\{ y \leq \mathcal{M} \right\} \right)
        > \frac{\epsilon}{C n^\nu} \right\}
    \leq 2 \exp \left( - \frac{2c}{C^\mu} n^{1 - \mu \nu} \epsilon^\mu \right) \, \text{;}
\end{multline}
\begin{setlength}{\multlinegap}{0pt}
    \begin{multline}\label{uc2}
\mathscr{V}_K \left\{ |x| \geq \mathcal{M} + \epsilon \right\}
    \leq \mathscr{V}_K \left\{ d_2 \left( x, \left\{ y \geq \mathcal{M} \right\} \right)
        > \epsilon \right\} \\
    \leq \mathscr{V}_K \left\{ d_K \left( x, \left\{ y \geq \mathcal{M} \right\} \right)
        > \frac{\epsilon}{C n^\nu} \right\}
    \leq 2 \exp \left( - \frac{2c}{C^\mu} n^{1 - \mu \nu} \epsilon^\mu \right) \, \text{;}
    \end{multline}
\end{setlength}

conclude with Theorems \ref{thm_aver}, \ref{intthm} and \ref{locthm} as in the proof of Corollary \ref{lp}.
\end{proof}

\appendix

\section{Proofs of technical lemmata}\label{S::sphproofs}

Here we prove Lemmata \ref{sph}, \ref{sphder} and \ref{logder}; the proofs are also
rather technical.

\vspace{10pt}

\begin{proof}[Proof of Lemma \ref{sph}]
First,
\begin{equation*}\begin{split}
\lim_{n \to \infty} \frac{\Gamma(n/2)}{\sqrt{\pi n} \, \Gamma((n-1)/2)}
    = \lim \frac{(n/2e)^{n/2}}{\sqrt{\pi n} \, ((n-1)/2)^{(n-1)/2}} \\
    = \lim \left(\frac{n}{n-1}\right)^{(n-1)/2} \times \sqrt{\frac{n}{2e\pi n}} = (2\pi)^{-1}
\end{split}\end{equation*}
by Stirling's  formula.

Now,
$$\exp(-\epsilon-\epsilon^2/2) \leq 1 - \epsilon \leq \exp(-\epsilon-\epsilon^2/(2(1-\epsilon)^2)) \, \text{;}$$
hence
\begin{multline*}
(1 - t^2/n)^{(n-3)/2} \, e^{t^2/2}
= \left\{ (1 - t^2/n) \, e^{t^2/(n-3)} \right\}^{(n-3)/2} \\
\leq \left\{ \exp(- t^2/n - t^4/(2n^2(1 - t^2/n)^2)) \: \exp(t^2/(n-3)) \right\}^{(n-3)/2} \\
\leq \exp \left( \frac{3}{2n} \, t^2 - \frac{n-3}{16n^2} \, t^4 \right)
\leq \exp \left( \frac{3}{2n} \, t^2 - \frac{1}{64n} \, t^4 \right)
\end{multline*}
for $n \geq 4$. For $t \leq 16$ and $n$ large enough $3t^2/2n \leq (1 + \epsilon)$; for $t \geq 16$
we have $3t^2/2n \leq 3t^4/256n$ and hence $3t^2/2n - t^4/64n < - t^4/256n$.

Similarly,
\begin{multline*}
(1 - t^2/n)^{(n-3)/2} \, e^{t^2/2}
= \left\{ (1 - t^2/n) \, e^{t^2/(n-3)} \right\}^{(n-3)/2} \\
\geq \left\{ \exp(- t^2/n - t^4/(2n^2)) \: \exp(t^2/(n-3)) \right\}^{(n-3)/2}
\geq \exp(-t^4/4n) \, \text{.}
\end{multline*}

This proves the first pair of inequalities; thereby
\begin{multline*}
1 - \Psi_n(t) \: = \: \int_t^\infty \psi_n(u) \, du
\: \leq \: (1 + \epsilon_n) \int_t^\infty \phi(u) \, \exp(-u^4/256n) \, du \\
\leq (1 + \epsilon_n) \, \exp(-t^4/256n) \int_t^\infty \phi(u) \, du \\
\: \leq \: (1 + \epsilon_n) \, (1 - \Phi(t) \, \exp(-t^4/256n) \, \text{.}
\end{multline*}

Similarly,
\begin{multline*}
1 - \Psi_n(t) - (1 - \epsilon_n) \,(1 - \Phi(t)) \, \exp(-t^4/324 n) \\
\geq (1 - \epsilon_n) \, \int_t^\infty \phi(u) \, \left[\exp(-u^4/4n) - \exp(-t^4/324n) \right] \, du \\
= (1 - \epsilon_n) \, \left[ \int_t^{2t} + \int_{2t}^{3t} + \int_{3t}^\infty\right]  \, \text{.}
\end{multline*}
The second integral is positive; integrating by parts, we see that the third integral equals
\begin{multline*}
- \int_{3t}^\infty (1 - \Phi(u)) \, \exp(- u^4/4n) \, \frac{u^3}{n} \, du \\
\geq - \exp(-t^4/324n) \int_{3t}^\infty (1 - \Phi(u)) \, \frac{u^3}{n} \, du \, \text{.}
\end{multline*}
The first one is at least
\begin{multline*}
\left[ \exp(-t^4/64n) - \exp(-t^4/324n) \right] \, \int_t^{2t} \phi(u) \, du \\
\geq \left[ \exp(-t^4/64n) - \exp(-t^4/324n) \right] \, \int_t^{2t} (1 - \Phi(u)) \, u \, du \\
\geq \left[ \exp(-t^4/64n) - \exp(-t^4/324n) \right] \, \int_t^{2t} (1 - \Phi(u)) \, \frac{u^3}{n} \, du
\end{multline*}
for $t^2 < n/4$. If $t \geq t_0$ this proves the remaining inequality ($t_0$
does not depend on $n$). $1 - \Psi_n(t) \rightrightarrows 1 - \Phi(t)$ on
$[0, \, t_0]$ and for these $t$ one can ignore $\exp(-t^4/n)$ in all the
expressions; hence the inequality also holds.
\end{proof}

Now we prove Lemma \ref{sphder}. The proof uses Lemma \ref{logder} that is proved
further on (without using Lemma \ref{sphder}).

\begin{proof}[Proof of Lemma \ref{sphder}]
By definition,
$$ \frac{1 - \Psi_n(t)}{\psi_n(t)} = \int_t^\infty \frac{\psi_n(s)}{\psi_n(t)} \, ds
    = t \int_0^\infty \frac{\psi_n((1+u)t)}{\psi_n(t)} \, du \, \text{.} $$

To obtain the upper bound, just note that if $2t^2 < n$, then by Lemma \ref{logder}
$$ \psi_n((1+u)t)/\psi_n(t) \leq \exp ( - C_1 \, u \, t^2 ) $$
and hence the integral is bounded by $(C_1t^2)^{-1}$.

For the lower bound restrict the integral to $[0, \, 1]$; if $8 t^2 < n$, $2 (1 + u)^2 t^2 < n$ and
the subintegral expression is bounded from below by $\exp(-C_2ut^2)$ on this interval. Hence the integral is not
less than $\frac{1 - \exp(-C_2t^2)}{C_2 t^2}$. This concludes the proof for $t > t_0$ (for a constant $t_0$
independent of $n$); for $0 < t_0$ use Gaussian approximation for $\psi_n$ and $\Psi_n$ (Lemma \ref{sph}) to verify
the inequality.

The second statement can be verified by formal differentiation.
\end{proof}

\begin{proof}[Proof of Lemma \ref{logder}]
\begin{multline*}
\psi_n(t) \, /\, \psi_n \bigl((1+u)\, t\bigr)
    = \left(\frac{1 - t^2/n}{1 - (1+u)^2 t^2/n}\right)^{(n-3)/2} \\
    = \left(1 + \frac{t^2 \, (2u + u^2)}{n - (1+u)^2 \, t^2} \right)^{(n-3)/2}
    \leq \exp \left(\frac{(u + u^2/2) \, t^2}{1-(1+u)^2t^2/n}\right) \\
    \leq \exp(3 u t^2) \quad \text{for $\frac{(1+u)^2\,t^2}{n} < 1/2$}
\end{multline*}
and hence the second inequality holds. On the other hand,
\begin{multline*}
\psi_n \bigl((1+u)\,t\bigr)\,/\,\psi_n(t)
    = \left(\frac{1 - t^2 \, \frac{(1+u)^2}{n}}{1 - t^2/n} \right) \\
    = \left(1 - \frac{t^2}{n-t^2} \, (2u + u^2) \right)^{(n-3)/2}
    \leq \exp \left( - (u + u^2/2) \, t^2 \, \frac{1-3/n}{1-t^2/n}\, \right) \\
    \leq \exp(- 6 u t^2) \quad \text{for $\frac{t^2}{n} < 1/2$ and $n > 6$}
\end{multline*}
and therefore the first inequality holds as well.
\end{proof}

\end{document}